\documentclass[11pt]{article}

\usepackage{latexsym}
\usepackage{amssymb}

\usepackage{amsopn}

\newcommand{\St}{\mathop {\rm St}\nolimits} 
\newcommand{\Int}{\mathop{\rm Int}\nolimits}

\newcommand {\IZ}{\mathbb{Z}}
\newcommand {\IQ}{\mathbb{Q}}
\newcommand {\IP}{\mathbb{P}}
\newcommand {\IR}{\mathbb{R}}

\newcommand {\qed} {$\Box$}

\DeclareMathSymbol{\Upset}{\mathop}{symbols}{"2A}
\DeclareMathSymbol{\Downset}{\mathop}{symbols}{"2B}
\DeclareMathSymbol{\upset}{\mathop}{symbols}{"22}
\DeclareMathSymbol{\downset}{\mathop}{symbols}{"23}

\newtheorem{theorem}{Theorem}[section]

\newtheorem{proposition}[theorem]{Proposition}
\newtheorem{corollary}[theorem]{Corollary}
\newtheorem{lemma}[theorem]{Lemma}
\newtheorem{example}[theorem]{Example}

\newtheorem{remark}[theorem]{Remark}
\newtheorem{question}[theorem]{Question}

\begin{document}

\begin{center}

{\large{\bf{A Note on Monotonically Metacompact Spaces}}}

\medskip by

\bigskip
Harold R. Bennett{\footnote{Texas Tech University, Lubbock,
    TX 79409}}, Klaas Pieter Hart{\footnote{TU Delft, Delft, Netherlands}},
and David J. Lutzer{\footnote{College of William and Mary,
    Williamsburg, VA 23187}}

\end{center}

\medskip
\noindent {\bf{Abstract}}:  We show that any metacompact Moore space is
monotonically metacompact and use that result to characterize monotone
metacompactness in certain generalized ordered (GO)spaces.  We show, for
example, that a generalized ordered space with a
$\sigma$-closed-discrete dense subset is metrizable if and only if it
is monotonically (countably) metacompact, that a monotonically
(countably) metacompact GO-space is hereditarily paracompact, and that
a locally countably compact GO-space is metrizable if and only if it
is monotonically (countably) metacompact. We give an example of a
non-metrizable LOTS that is monotonically metacompact, thereby
answering a question posed by S. G.  Popvassilev. We also give
consistent examples showing that if there is a Souslin line, then
there is one Souslin line that is monotonically countable metacompact,
and another Souslin line that is not monotonically countably
metacompact.

\medskip
\noindent
{\bf{Keywords}}: metacompact, countably metacompact, monotonically
metacompact, monotonically countably metacompact, generalized ordered space,
GO-space, LOTS, metacompact Moore space, metrizable,
$\sigma$-closed-discrete dense set.

\medskip
\noindent
{\bf{MR numbers}}: Primary = 54D20; Secondary = 54E30, 54E35, 54F05

\bigskip
******************draft of October 6, 2009

\medskip

\section{Introduction}

For two collections of sets $\cal U$ and $\cal V$
we write ${\cal U} \prec {\cal V}$ to mean that for
each $U \in {\cal U}$ there is some $V \in {\cal
  V}$ with $U \subseteq V$.  Clearly ${\cal U} \prec {\cal V}$ implies
$\bigcup {\cal U} \subseteq \bigcup {\cal V}$ but it might happen
that $\bigcup {\cal U} \not= \bigcup {\cal V}$.

\medskip A space $X$ is {\em{(countably) metacompact}} if each
(countable) open cover of $X$ has a point-finite open refinement that
also covers $X$.  Popvassilev \cite {Popvassilev} defined that a space
is {\em{ monotonically (countably) metacompact}} if there is a
function $r$ that associates with each (countable) open cover $\cal U$
of $X$ an open point-finite refinement $r({\cal U})$ that covers $X$,
where $r$ has the property that if ${\cal U}$ and $\cal V$ are open
covers with ${\cal U} \prec {\cal V}$ then $r({\cal U}) \prec r({\cal
  V})$.  The function $r$ is called a {\em{monotone (countable)
    metacompactness operator}} for $X$. {\em{Warning: The literature
    contains other, non-equivalent definitions of monotone countable
    metacompactness{\footnote{The referee of this paper pointed out
        that Good, Knight, and Stares \cite {GoodKnightStarnes} have
        given another definition of a monotone countable
        metacompactness property.  According to Ying and Good
        \cite{YingGood}, that property is equivalent to Hodel's
        $\beta$-space property \cite{Hodel}.  Consequently, the
        Good-Knight-Stares and Popvassilev definitions of monotone
        countable metacompactness are {\em{not}} equivalent, and
        neither implies the other. On the one hand, the space
        $[0,\omega_1)$ with is usual topology is certainly a
        $\beta$-space but in the light of Proposition
        \ref{HereditaryParacompact} of this paper, $[0, \omega_1)$ is not
        monotonically countably metacompact in Popvassilev's sense. On
        the other hand, the example machine $Bush(S,T)$ in
        \cite{BennettLutzerBush} constructs a family of linearly
        ordered spaces that are always monotonically metacompact in
        the sense of Popvassilev, but never $\beta$-spaces.}} and
    throughout this paper we study the monotone metacompactness
    property of Popvassilev.}}

\medskip 

\noindent 

\medskip
Our first main result shows that monotone metacompactness is a
property of many generalized metric spaces. As a corollary to a more
technical Proposition \ref{MetaMoore} we will show:

\begin{theorem} \label{MetaMooreTheorem} Any metrizable space, and any
  metacompact Moore space, is monotonically metacompact. \end{theorem}

Recall that a {\em{GO-space}} is a triple $(X, \tau, <)$
where $<$ is a linear ordering of the set $X$ and $\tau$ is a
Hausdorff topology on $X$ that has a base of convex subsets of $(X,
<)$, possibly including some singletons.  For any GO-space $(X, \tau,
<)$, we have $\lambda \subseteq \tau$, where $\lambda$ is the
usual open interval topology of the ordering $<$. If $\tau =
\lambda$, then $(X, \lambda, <)$ is called a {\em{linearly ordered
    topological space}} (LOTS) and for any GO-space $(X, \tau,<)$, $(X,
\lambda,<)$ is called the {\em{underlying LOTS for the GO-space}}.

\medskip
 Our later results will deal with
GO-spaces $(X, \tau, <)$ that have dense sets that are
$\sigma$-closed-discrete, and with GO-spaces $(X, \tau, <)$ whose
underlying LOTS $(X, \lambda, <)$ has a dense subset that is
$\sigma$-closed-discrete in $(X, \lambda)$.  Of particular interest
are separable GO-spaces and GO-spaces whose underlying
LOTS is separable. The best-known examples of this type are the
Sorgenfrey and Michael lines, and GO-spaces constructed from the
Alexandroff double arrow, i.e., the set $X = \IR \times \{0,1\}$ with
the lexicographic ordering. 

\medskip 

B.J Ball \cite{Ball} has shown that any LOTS is
{\em{countably paracompact}} (= any countable open
cover has a locally finite refinement), and hence
countably metacompact, and that result was
extended to GO-spaces in
\cite{LutzerThesis}. Adding ``monotonicity'' to
the countable metacompactness property is a
significant strengthening, as our results will
show. For example, in Propositions \ref{HereditaryParacompact} and
\ref{CompactMonoMeta}, we prove:

\begin{theorem} \label{IntroHeredPara} Suppose $(X, \tau, <)$ is a
  GO-space that is monotonically countably metacompact.  Then $(X,
  \tau)$ is hereditarily paracompact. \end{theorem}

\begin{theorem} \label{IntroCompact} Suppose $(X, \tau, <)$ is a
  compact LOTS or, more generally, a locally countably compact
  GO-space.  If $(X, \tau)$ is
  monotonically countably metacompact, then $X$ is
  metrizable. \end{theorem}

\medskip 
In certain cases we can characterize which
GO-spaces are monotonically (countably)
metacompact. To state our result, we need some
special notation. For any GO-space $(X, \tau, <)$,
let

\begin{quote} 

$I_\tau := \{x \in X: \{x\} \in \tau  \}$,

$R_\tau := \{ x \in X - I_\tau : [x \rightarrow) \in \tau \}$,

$L_\tau := \{x \in X - I_\tau : (\leftarrow, x] \in \tau\}$, and 

$E_\tau := X - (I_\tau \cup R_\tau \cup L_\tau)$.

\end{quote}

\noindent {\underbar{Warning}}: The above definitions for sets
$I_\tau, R_\tau, L_\tau$, and $E_\tau$ are slightly different from
definitions given in other papers on GO-spaces. In some other papers, the
set of right-looking points is defined as $\{ x \in X: [x,\rightarrow)
\in \tau - \lambda\}$, but theorems in this paper require that
$R_\tau$ must be defined as above.

\medskip

\begin{theorem} \label{MainTheorem} Let $(X, \tau, <)$ be a GO-space
  whose underlying LOTS $(X, \lambda, <)$ has a
  $\sigma$-closed-discrete dense subset.  Then the following are
  equivalent: 

\begin{itemize}

\item [{a)}] $(X, \tau)$ is monotonically
  metacompact; 

\item [{b)}] $(X, \tau)$ is monotonically countably metacompact;

\item [{c)}] the set $R_\tau \cup L_\tau$ is
  $\sigma$-closed-discrete in $(X, \tau)$;

\item [{d)}] the set $R_\tau \cup L_\tau$ is
  $\sigma$-closed-discrete in $(X, \lambda)$. 

\end{itemize} 

\end{theorem}

\noindent
Theorem \ref{MainTheorem} applies to GO-spaces constructed on the
usual real line and shows that the Michael line is monotonically
countably metacompact, while the Sorgenfrey line is not. In addition,
it shows that the Alexandroff double arrow is not monotonically
metacompact.

\medskip

The proof of Theorem \ref{MainTheorem}, combined with M.J. Faber's
metrization theorem for GO-spaces \cite{Faber}, will show:

\begin{corollary} \label{Metrization} Let $(X, \tau, <)$ be a GO-space
  with a $\sigma$-closed-discrete dense subset.  Then the following
  are equivalent:

\begin{itemize}

\item [{a)}] $(X, \tau)$ is
  monotonically metacompact; 

\item [{b)}] $(X, \tau)$ is monotonically countably metacompact;

\item [{c)}] $(X, \tau)$ is metrizable.

\end{itemize}

\end{corollary}

\noindent
Our next example, based on
the Michael line (see Example \ref{MonotmetaLOTS} for details),
shows that the hypothesis in Corollary
\ref{Metrization} concerning the existence of a
$\sigma$-closed-discrete dense subset cannot be
removed, and at the same time, it answers
a question posed in a recent paper \cite{ Popvassilev} by
S. G. Popvassilev, namely
``Must a monotonically metacompact LOTS be
metrizable?''

\begin{example} \label{MichaelExample} There is a LOTS that is
  monotonically metacompact but not metrizable. \end{example}

 Theorem
\ref{MainTheorem} can also be used to show that in certain GO-spaces,
monotonic (countable) metacompactness is a hereditary property.

\begin{corollary} \label{AlmostHereditary} Suppose $(X, \tau, <)$ is a
  GO-space whose underlying LOTS $(X, \lambda, <)$ has a
  $\sigma$-closed-discrete dense set. Let $Y \subseteq X$. Then the
  subspace $(Y, \tau_Y, <)$ is monotonically (countably) metacompact.
\end{corollary}

\noindent
However, we do not know whether monotone (countable) metacompactness
is a hereditary property in other kinds of GO-spaces. 

\medskip

To what extent can our results be extended?  It is
known that any GO-space with a
$\sigma$-closed-discrete dense subset is
{\em{perfect}} (= every closed set is a
$G_\delta$-set in the space) and one might ask,
for example, whether Corollary \ref{Metrization}
could be proved for perfect GO-spaces. We will
show that the answer is ``Consistently no'' by
looking at a {\em{Souslin line}} (= a LOTS that
satisfies the countable chain condition but is not
separable). Souslin lines are hereditarily
Lindel\"{o}f and therefore perfect. Whether such
spaces exist is undecidable in ZFC
\cite{Todorcevic}. Whether there is a perfect LOTS
that does not have a $\sigma$-closed-discrete
dense subset is an old problem of Maarten Maurice
that is undecidable in ZFC, at least for spaces of
weight $\omega_1$, and is intimately related to
the Souslin problem
\cite{BennettLutzerMaurice}. In our paper's final
section, we show:

\begin {example} \label{SouslinExample} If there is a Souslin line, then
  some Souslin lines are monotonically countably
  metacompact, while other Souslin lines are not. \end{example}

Throughout this paper,  $\IR,~ \IP$, and $ \IQ$ denote the sets of real,
irrational, rational numbers respectively, and $\IZ$ is the set of
integers. The authors want to thank the referee for a series of
helpful remarks that improved the current paper and suggested
directions for further investigation (see Question \ref{RefereeQuestion}).

\section{Preliminary results} 

\medskip

We must carefully distinguish between subsets of a space that are
{\em{relatively discrete}} (= their subspace topology is the discrete
topology) and subsets that are {\em{closed-discrete}} (= every point of
the space has a neighborhood containing at most one point of the given
subset). Clearly, a set is closed-discrete if it is both relatively
discrete and closed. In general, we will need to distinguish between
subsets that are {\em{$\sigma$-relatively-discrete}} (= countable
unions of relatively discrete subsets) and those that are
{\em{$\sigma$-closed-discrete}} (= countable unions of closed-discrete
subspaces). However, the two concepts are equivalent in perfect
    spaces as our next lemma
shows. The lemma is well-known and easily proved, and applies to any
topological space, not just to GO-spaces.

\begin{lemma} \label{Perfect} If $(X,\tau)$ is a perfect topological
  space, then any relatively discrete subset is
  $\sigma$-closed-discrete. Hence, any $\sigma$-relatively discrete
  subset of a perfect space $(X, \tau)$ is $\sigma$-closed-discrete in
  $(X, \tau)$.
  \end{lemma}

\noindent
Proof: The proof is a standard argument but we include it for completeness.
Let $D$ be a relatively discrete subset of $X$ and for each $x \in D$
let $U(x)$ be an open set with $U(x) \cap D = \{x\}$. Write the open
set $V := \bigcup \{U(x): x \in D \}$ as $V = \bigcup \{ F(n): n \geq
1 \}$ where each $F(n)$ is closed in $X$.  Then the set $D(n):= D \cap
F(n)$ is closed and discrete and $D = \bigcup \{ D(n): n \geq 1 \}$. \qed

\begin{lemma} The existence of a $\sigma$-closed-discrete dense
set in a GO-space $(X, \tau, <)$ is a hereditary property 
 and implies that $(X, \tau)$
is perfect. \cite{BennettLutzerPurisch}   \qed \end{lemma}

By way of contrast, the existence of a
$\sigma$-relatively discrete dense set in a
GO-space $(X, \tau,<)$ is not enough to make $(X,
\tau)$ perfect and is not a hereditary
property. For example, in the Michael line $M$,
the set of irrationals is a relatively discrete
dense set, but $M$ is not perfect. Example 5.3 in 
\cite{BennettLutzerPurisch} describes a GO-space that has a
$\sigma$-relatively-discrete dense set, and has a subspace that does
not.

\begin{lemma} \label{MakingDiscrete} Let $E$ be a
  closed-discrete subset of a GO-space $(X, \tau, <)$ and let $S
  \subseteq X$.  Suppose that for each $x \in S$ there is some $e(x)
  \in E$ with $x < e(x)$ and such that the collection ${\cal C} := \{
  [x, e(x)] : x \in S \}$ is a pairwise disjoint collection.  Then the
  collection ${\cal C}$ is discrete in $(X, \tau)$ and the set $S$ is
  a closed-discrete subset of $(X, \tau)$. \end{lemma}

\noindent
Proof: Let $y \in X$ and let $U$ be a convex neighborhood of $y$ that
contains at most one point of $E$. Suppose $[x_i, e(x_i)]$ are four
distinct members of ${\cal C}$ and for contradiction suppose that $U$
meets all four sets.  Without loss of generality, we may assume $x_1 <
x_2 < x_3 < x_4$. Because the collection $\cal C$ is pairwise
disjoint, we must have $x_1 \leq e(x_1) < x_2 \leq e(x_2) < x_3 \leq
e(x_3) < x_4 \leq e(x_4)$. Then convexity of $U$, plus that fact that
$U$ meets both $[x_1, e(x_1)]$ and $[x_4, e(x_4)]$, shows that both
$e(x_2)$ and $e(x_3)$ belong to $U$ and to $E$, and that is
impossible.  Therefore the open set $U$ meets at most three members of
$\cal C$.  Because $\cal C$ is pairwise disjoint, that is enough to
show that $\cal C$ is a discrete collection. Because $S$ contains
exactly one point from each member of a discrete
collection, the set $S$ is closed and discrete. \qed

\begin{lemma} \label{SigmaDiscrete} Suppose $(X, \tau, <)$ is a
  GO-space and that the underlying LOTS $(X, \lambda)$ has a
  $\sigma$-closed-discrete dense set $D = \bigcup \{( D(n): n \geq 1
  \}$.  Suppose $S \subseteq X$ is $\sigma$-relatively discrete in the
  subspace topology $\tau_S$ and that no point of $S$ is isolated in
  $\tau$. Then $S$ is $\sigma$-closed-discrete in $(X, \lambda)$  and
  therefore also in $(X, \tau)$. \end{lemma}

\noindent
Proof: It is enough to prove the lemma in case $S$ is relatively
discrete in $(X, \tau)$. For each $x \in S$ let $U(x)$ be a convex
$\tau$-neighborhood of $x$ with $U(x) \cap S = \{x\}$. Because $S \cap
I_\tau = \emptyset$, we have $S \subseteq R_\tau \cup E_\tau \cup
L_\tau$. Let $S_1 := S \cap (R_\tau \cup E_\tau)$. For each $x \in
S_1$ there is some $y(x)> x$ with $[x, y(x)) \subseteq U(x)$.  Because
$x \in R_\tau \cup E_\tau$, the set $(x, y(x))$ is not empty. Then
there is some integer $N(x)$ and some $d(x) \in D(N(x)) \cap (x,
y(x))$. Replacing $y(x)$ by some point of $(x, y(x))$ if necessary, we
may assume that $d(x)$ is the only point of $(x, y(x)) \cap D(N(x))$.
By Lemma \ref{MakingDiscrete}, for each $k \geq 1$ the collection $\{
[x, d(x)] : N(x) = k \}$ is a discrete collection in $(X, \lambda)$,
and $S_1(k) := \{x \in S_1: N(x) = k \}$ is a closed-discrete subset
of $(X, \lambda)$ for each $k$. But $S_1 = \bigcup \{S_1(k): k \geq 1
\}$.  Points of $S \cap L_\tau$ are treated analogously. Hence $S$ is
$\sigma$-closed-discrete in $(X, \lambda)$ and hence also in $(X,
\tau)$. \qed


\medskip
\section{Main theorems}

\begin{proposition} \label{MetaMoore} Suppose $X$ is a metacompact
  Moore space.  There is a function $r$ such that for each collection
  ${\cal U}$ of open subsets of $X$, $r({\cal U})$ is a collection
  of open subsets of $X$ satisfying:

\begin{itemize}

\item [{a)}] $r({\cal U})$ is point-finite;

\item [{b)}] $r({\cal U}) \prec {\cal U}$;

\item [{c)}] $\bigcup r({\cal U}) = \bigcup {\cal U}$;

\item [{d)}] if $G, H \in r({\cal U})$ have $G \subseteq H$ then $G =
  H$; and 

\item [{e)}] if ${\cal V}$ is an open collection
  with ${\cal U} \prec {\cal V}$, then $r({\cal U})
  \prec r({\cal V})$.

\end{itemize}

 \end{proposition}

\noindent
Proof:  Suppose $X$ is a metacompact Moore space.
Let $\langle {\cal G}(n) \rangle$ be a development for $X$
where each ${\cal G}(n)$ is a point-finite open
cover of $X$.  We may assume that ${\cal G}(n+1)
\prec {\cal G}(n)$.  In addition, for each $n$, every point of
$X$ belongs to a maximal member of the point-finite collection ${\cal
  G}(n)$  so we may
assume that each member of each ${\cal G}(n)$ is maximal in ${\cal
  G}(n)$, i.e., if $G, H \in {\cal G}(n)$ are distinct, then neither $G
\subseteq H$ nor $H \subseteq G$.

\medskip Let $\cal U$ be a collection of open subsets of $X$. Define
${\cal U}(1) = \{G \in {\cal G}(1): G \prec {\cal U}\}$, where we
write $G \prec {\cal U}$ to mean that $G$ is a subset of some member
of $\cal U$. For $n \geq 1$, let $${\cal U}(n+1) := \{ G \in {\cal
  G}(n+1): G \prec {\cal U} ~{\rm{and}}~ G \not\prec \bigcup \{ {\cal
  U}(i): 1 \leq i \leq n \} \}.$$
Then $r({\cal U}) := \bigcup \{
{\cal U}(n): n \geq 1 \}$ is a collection of open sets in $X$ that
refines $\cal U$ and has $\bigcup r({\cal U}) =
\bigcup {\cal U}$. 

\medskip
Next we show that $r({\cal U})$ is point-finite. Fix $p \in X$
with $p \in \bigcup {\cal U}$. Find the first $n$ such that $p \in
\bigcup {\cal U}(n)$.  Then we have some $G_n \in {\cal G}(n)$ with $p
\in G_n$ where $G_n \prec {\cal U}$. Find $m \geq n+1$ so that $\St(p,
{\cal G}(m)) \subseteq G_n$ and note that if $k \geq m$ and $p \in G
\in {\cal G}(k)$, then $G \subseteq \St(p, {\cal G}(m)) \subseteq G_n$
so that $G \prec \bigcup \{ {\cal U}(j) : j \leq m \} \}$ and
therefore $G \not\in {\cal U}(k)$. Consequently, $$\{G \in r({\cal
  U}): p \in G \} \subseteq \bigcup \{ {\cal G}(j): j \leq m \}$$
and the latter collection is point-finite.  Hence $r({\cal U})$ is
also point-finite. 

\medskip To prove (d), suppose distinct $G, H \in r({\cal U})$ have $G
\subset H$. Find integers $m, n$ with $G \in {\cal G}(m)$ and $H \in
{\cal G}(n)$. Then $m \not= n$ because each member of ${\cal G}(n)$ is
maximal. We cannot have $m > n$ because no member of ${\cal G}(m)$ was
chosen for $r({\cal U})$ if it was contained in a previously chosen
member of $r({\cal U})$.  So we must have $m < n$.  Because $H \in
{\cal G}(n) \prec {\cal G}(m)$ there is some $G' \in {\cal G}(m)$ with
$H \subseteq G'$.  But then $G \subset H \subseteq G'$, which shows that the
element $G \in {\cal G}(m)$ is not maximal in ${\cal G}(m)$, and that is
impossible. Hence (d) holds.

\medskip
To verify (e), suppose ${\cal U} \prec {\cal V}$.  Clearly ${\cal
  U}(1) \prec {\cal V}(1)$. Suppose $n \geq 1$ and that $\bigcup \{
{\cal U}(i) : i \leq n \} \prec \bigcup \{ {\cal V}(i) : i \leq n
\}$. Let $G \in {\cal U}(n+1)$.  Then $G \prec {\cal U}$ so that $G
\prec {\cal V}$.  If $G \prec \bigcup \{ {\cal V}(i): i \leq n \}$
then $G \prec r({\cal V})$, and otherwise $G \in r({\cal V})$. Hence
(e) holds.   \qed

\medskip An immediate corollary of the previous result is Theorem \ref
{MetaMooreTheorem}.

\medskip
\begin{corollary} \label{MonotMetaCorollary}  Any
  metacompact Moore space, and any metric space,
  is monotonically metacompact.
  \qed \end{corollary}

\begin{corollary} \label{DiscretizingMetric} Suppose $(X, \mu)$ is a
  metrizable or metacompact Moore space and $S \subseteq X$.  Let
  $\mu^S$ be the topology on $X$ having $\mu \cup \{ \{x\}: x \in S
  \}$ as a base. Then $(X, \mu^S)$ is monotonically metacompact.
\end{corollary}

\noindent
Proof: By Proposition \ref{MetaMoore} we know that the 
space $(X, \mu)$ has a monotone metacompactness operator $r$ that acts
on collections of $\mu$-open sets, even if they do not cover $X$. Let
${\cal U}$ be any open cover of $(X, \mu^S)$. Define ${\cal U}_\mu :=
\{ \Int_\mu(U): U \in {\cal U} \}$ and note that $X - S \subseteq
\bigcup {\cal U}_\mu$. Find the point-finite $\mu$-open refinement
$r({\cal U}_\mu)$ and define $s({\cal U}) := r({\cal U}_\mu) \cup \{
\{x\}: x \in S \}$.  Then the collection $s({\cal U})$ is point-finite
in $X$, covers all of $X$, and refines ${\cal U}$.  Further, if ${\cal
  U} \prec {\cal V}$ then $s({\cal U}) \prec s({\cal V})$ as
required. \qed

\medskip

Experience has shown that adding ``monotonicity''
to a covering property makes the property much
stronger.  The best example of this is Gary Gruenhage's
proof that a monotonically compact{\footnote{$X$
    is {\em{monotonically compact}} if for every
    open cover ${\cal U}$ of $X$, there is a finite
    open refinement $r({\cal U})$ such that if
    ${\cal U} \prec {\cal V}$ then $r({\cal U})
    \prec r({\cal V})$. }} Hausdorff space must be
metrizable \cite{Gruenhage}.  As noted in the
Introduction, every GO-space is countably
metacompact.  Our next result (which is Theorem \ref{IntroHeredPara} of
the Introduction) shows that adding
monotonicity to countable metacompactness 
makes the property much stronger.

\begin{proposition} \label{HereditaryParacompact} Suppose $(X, \tau,
  <)$ is a GO-space that is monotonically countably metacompact. Then
  $(X, \tau)$ is hereditarily paracompact. \end{proposition}

\noindent
Proof: If $(X, \tau)$ is not hereditarily paracompact, then by
\cite{EngelkingLutzer} there is an uncountable regular cardinal
$\kappa$ and a stationary subset $S \subseteq [0, \kappa)$ that embeds
in $X$ under a mapping that is strictly increasing, or strictly
decreasing.  Consider the case where the mapping is strictly
increasing, the other case being analogous. Then we may view $S$ as a
subset of $X$ and know that the ordering $<_S$ inherited from $(X,<)$
is the same as the ordering of $S$ as a subspace of $\kappa$. This
allows us to write such things as ``if $\alpha \in S$, then $\alpha^+
\in X$,'' where $\alpha^+$ is the first element of $S$ that
lies above $\alpha$, and ``$\{ (\leftarrow, \alpha^+ ), (\alpha,
\rightarrow)\}$ is a open cover of $X$''.

\medskip Suppose there is a monotone countable metacompactness
operator $r$ on $X$.  Let $S^d$ be the set of limit points of $S$ in
$X$ that belong to $S$.  Then $S^d$ is also stationary in $\kappa$.
For each $\alpha \in S^d$ consider the open cover ${\cal U}(\alpha) =
\{(\leftarrow,\alpha^+), (\alpha, \rightarrow)\}$ of $X$ and find $r({\cal
  U}(\alpha))$. Choose $O(\alpha) \in r({\cal U}(\alpha))$ with
$\alpha \in O(\alpha)$. Then $O(\alpha) \subseteq (\leftarrow,
\alpha^+)$ and there is some $f(\alpha) \in S$ with $f(\alpha) <
\alpha$ such that $[f(\alpha), \alpha] \subseteq O(\alpha)$.  The
Pressing Down Lemma provides some $\beta \in S$ such that the set $T
:= \{\alpha \in S^d: f(\alpha) = \beta \}$ is stationary.  Choose a
strictly increasing sequence $\alpha(1) < \alpha(2) < \cdots$ in $T$
with the property that $\alpha(n)^+ < \alpha(n+1)$ and let ${\cal V} =
\bigcup \{ {\cal U}(\alpha(n)): n \geq 1 \}$. Then ${\cal V}$ is a
countable open cover of $X$ so that $r({\cal V})$ is defined. For each
$i$, the cover $r({\cal U}(\alpha(i)))$ refines $r({\cal V})$ so there is
some $W(i) \in r({\cal V})$ with $O(\alpha(i)) \subseteq W(i)$. Note
that for each $i$ we have $$\beta \in [\beta, \alpha(i)]=
[f(\alpha(i)), \alpha(i)] \subseteq O(\alpha(i)) \subseteq W(i),$$
because $\alpha(i) \in T \subseteq S^d$ gives $f(\alpha(i)) = \beta$.
We will show that there are infinitely many distinct sets in the
collection $\{W(i): i \geq 1 \}$ and that will contradict
point-finiteness of $r({\cal V})$.

\medskip Let $j_1 = 1$ and consider $\alpha(j_1) \in O(\alpha(j_1))
\subseteq W(j_1)$.  Because $W(j_1) \in r({\cal V)}$ which refines
${\cal V} = \bigcup \{{\cal U}(\alpha(i)): i \geq 1 \}$ there is some
$i_1 \geq 1$ for which either $W(j_1) \subseteq (\leftarrow,
a(i_1)^+)$ or else $W(j_1) \subseteq (\alpha(i_1) , \rightarrow)$. The
second alternative cannot happen because $W(j_1)$ contains $\beta <
\alpha(i_1)$ while $(\alpha(i_1), \rightarrow)$ does not, so we have
$W(j_1) \subseteq (\leftarrow, \alpha(i_1)^+)$. Let $j_2 = i_1 +1$.
Note that $\alpha(i_1)^+ < \alpha(j_2)$ and consider $\alpha(j_2)
\in W(j_2)$.  Because $\alpha(j_2) \not\in W(j_1)$ we know that
$W(j_2) \not= W(j_1)$. Because $W(j_2) \in r({\cal V})$ and $r({\cal
  V})$ refines $\cal V$, there is some $i_2$ such that either $W(j_2)
\subseteq (\leftarrow, \alpha(i_2)^+)$ or else $W(j_2) \subseteq
(\alpha(i_2), \rightarrow)$.  The second alternative cannot occur
because $\beta \in [\beta, \alpha(j_2)] \subseteq O(\alpha(j_2))
\subseteq W(j_2)$ while $\beta \not\in (\alpha(i_2), \rightarrow)$, so
$W(j_2) \subseteq (\leftarrow, \alpha(i_2)^+)$. Let $j_3 = i_2 + 1$ and
consider $W(j_3)$.  Because $\alpha(j_3) \in W(j_3)$ and $\alpha(j_3)
\not\in W(j_1) \cup W(j_2)$ we see that the sets $W(j_1), W(j_2)$, and
$W(j_3)$ are distinct.  This recursion continues, producing an
infinite sequence of distinct members $W(j_k)$ of $r({\cal V})$, with $\beta
\in W(j_k)$ for each $k$, and that is impossible because $r({\cal V})$
is point-finite.  \qed

\medskip
There is a generalization of Proposition \ref{HereditaryParacompact}
that might be of interest. A deep result of Balogh and Rudin \cite{BaloghRudin}
  shows that a monotonically normal space is paracompact if and only
  if it does not contain a closed subspace that is a topological copy
  of a stationary subset of a regular uncountable cardinal.
  Therefore, because monotone countable metacompactness is a
  closed-hereditary property, the proof given for the previous theorem
  actually shows that a monotonically normal space that is
  monotonically countably metacompact must be paracompact.

\medskip

Popvassilev proved in \cite{Popvassilev} that neither of the ordinal
spaces $[0, \omega_1)$ and $[0, \omega_1]$ is monotonically countably
metacompact. Our Proposition \ref{HereditaryParacompact} gives another
proof of that result. 

\medskip
One might wonder whether, among subspaces of ordinals, the hypothesis
of monotone metacompactness would give a
conclusion even stronger than hereditary paracompactness. Our next
example shows that one cannot obtain metrizability from monotone 
metacompactness. Because our next example is a LOTS under a different
order, it solves Popvassilev's question from \cite {Popvassilev}, but
it fails to be first-countable.  A first-countable example is given in
Example \ref{MonotmetaLOTS}, below.

\begin{example} There is a subspace $X \subseteq [0, \omega_1]$ such
  that $X$ is monotonically metacompact but not metrizable. \end{example}

\noindent
Proof: Let $X := \{\alpha \in [0, \omega_1): \alpha~
{\rm{is~not~a~limit~ordinal}} \} \cup \{\omega_1\}$ topologized as a
subspace of $[0, \omega_1]$.  Then in its subspace topology, $X$ is a
GO-space. (Note that, under a different ordering, $X$ is actually a
LOTS.) Let $\cal U$ be any open cover of $X$.  Let $\beta$ be the
first ordinal such that $(\beta, \omega_1] \cap X$ is a subset of some
member of $\cal U$ and define $r({\cal U}) := \{ (\beta, \omega_1]
\cap X \} \cup \{ \{ \gamma\}: \gamma \leq \beta, \gamma \in X \}$.
Then $r$ is a monotone metacompactness operator for $X$, and yet $X$
is not metrizable.  \qed

\medskip
Our next result proves Theorem \ref{IntroCompact} of the Introduction.

\begin{proposition} \label{CompactMonoMeta} Suppose $(X, \lambda, <)$ is
  a compact LOTS. Then $X$ is monotonically
  countably metacompact if and only if $(X, \lambda)$ is
  metrizable.\end{proposition}

\noindent
Proof: If $X$ is metrizable, then $X$ is monotonically metacompact by
Corollary \ref{MonotMetaCorollary}. To prove the other half, suppose
$(X,\lambda)$ is compact and monotonically countably metacompact. We will
show that $X$ is {\em{monotonically countably compact}} (= every
countable open cover ${\cal U}$ has a finite open refinement $r({\cal U})$
in such a way that if ${\cal U}$ and ${\cal V}$ are countable open covers
with $\cal U$ refining $\cal V$ then $r({\cal U})$ refines $r({\cal
  V})$), and then we will invoke Popvassilev's theorem that any monotonically
countably compact LOTS is metrizable \cite{Popvassilev}.

\medskip

Let ${\cal U}$ be any countable open cover of $X$.
Monotone countable metacompactness gives a point-finite
refinement $r({\cal U})$ whose members are convex
subsets of $X$. Replace $r({\cal
  U})$ by the subcollection $s({\cal U}) := \{ V
\in r({\cal U}): V ~{\rm{is~maximal~in}}~ r({\cal
  U})\}$. If $r$ is a monotone metacompactness
operator, then so is $s$. Notice that no member of
$s({\cal U})$ is contained in any other member of
$s({\cal U})$.

\medskip We claim that $s({\cal U})$ is finite.  If not, choose
infinitely many distinct sets $V_i \in s({\cal U})$.  At most one
contains the left end-point of $X$, for otherwise one member of
$s({\cal U})$ would be contained in another member of $s({\cal U})$
and that is impossible, and at most one contains the right endpoint of
$X$. Discarding those two, we may assume that each $V_i = (a_i, b_i)$
for some $a_i, b_i \in X$. For each $i$, $V_i$ is the only member of
$s({\cal U})$ having $a_i$ as its left endpoint, because each member
of $s({\cal U})$ is maximal in $s({\cal U})$, so that $a_i \not= a_j$
whenever $i \not= j$.  Passing to a subsequence if necessary, we may
assume that $\langle a_i \rangle$ is a strictly monotone sequence.
Consider the case where $\langle a_i \rangle$ is increasing. Suppose
$i < j$.  If $b_j \leq b_i$ then we would have $a_i < a_j < b_j \leq
b_i$ and that is impossible because no member of $s({\cal U})$ can
contain another. Hence $\langle b_i \rangle$ is also monotone
increasing.  Because $X$ is compact, $p := \sup \langle a_i \rangle$
and $q := \sup \langle b_i \rangle$ are points of $X$ and $p \leq q$.
If $p < q$ then infinitely many members of the point-finite collection
$s({\cal U})$ contain $p$, so we have $p = q$.  Choose any $V \in
s({\cal U})$ that contains $p$.  Then $V$ contains infinitely many of
the sets $V_i = (a_i, b_i)$ and that is impossible because no member
of $s({\cal U})$ can contain another member of $s({\cal U})$. \qed

\begin{corollary} Suppose $(X, \tau, <)$ is
  a locally countably compact GO-space. Then $X$ is monotonically
  (countably) metacompact if and only if $(X, \tau)$ is
  metrizable. \end{corollary}

\noindent
Proof: If $(X, \tau)$ is metrizable, then it is monotonically
metacompact by Proposition \ref{MetaMoore} and therefore also
monotonically countably  metacompact. Next, suppose $(X, \tau, <)$ is
monotonically (countably) metacompact. By Proposition
\ref{HereditaryParacompact}, $(X, \tau)$ is hereditarily paracompact,
so that $(X, \tau)$ is locally compact.  Because of Proposition
\ref{CompactMonoMeta}, $(X, \tau)$ is locally metrizable.  Because
$(X, \tau)$ is paracompact, we see that $(X, \tau)$ is metrizable. \qed

\medskip Recall the special subsets $R_\tau$ and $L_\tau$ defined for
a GO-space $(X, \tau, <)$ in the Introduction. Theorems of M.J. Faber
\cite{Faber} and Jan van Wouwe \cite{vanWouwe} show that one key
to metrization theory for a GO-space $(X, \tau, <)$ is the
hypothesis that $R_\tau \cup L_\tau$ is a $\sigma$-closed-discrete
subset of $(X, \tau)$.  We will show that this same hypothesis plays a
central role in the study of monotone metacompactness in GO-spaces. We
will need the extra hypothesis that the underlying LOTS $(X, \lambda,
<)$ of the given GO-space has a $\sigma$-closed-discrete dense subset.
Examples in the final section of our paper will show that this extra
hypothesis is needed.

\medskip The remaining results in this section deal with GO-spaces
that have $\sigma$-closed-discrete dense subsets, and GO-spaces whose
underlying LOTS have $\sigma$-closed-discrete dense subsets. We will
combine them to prove Theorem \ref{MainTheorem} and Corollary
\ref{Metrization} of the Introduction. Upon first reading of these
results, it might be helpful for the reader to replace the hypothesis
``$\sigma$-closed-discrete dense set'' by ``separable''.

\begin{proposition} \label{NewFirstHalf} Suppose $(X, \tau, <)$ is a
  GO-space for which the underlying LOTS $(X, \lambda, <)$ has a
  $\sigma$-closed-discrete dense set. If $(X, \tau)$ is monotonically
  countably metacompact, then $R_\tau \cup L_\tau$ is
  $\sigma$-closed-discrete as a subspace of $(X, \tau)$ and as a
  subspace of $(X,\lambda)$. \end{proposition}

\noindent
Proof: From the definition of $R_\tau$ (in the Introduction), no point
of $R_\tau$ is isolated in $(X, \tau)$ so that Lemma
\ref{SigmaDiscrete} guarantees that $R_\tau$ is
$\sigma$-closed-discrete in $(X, \tau)$ if and only if $R_\tau$ is
$\sigma$-closed-discrete in $(X, \lambda)$.  The same assertion holds
for $L_\tau$.

\medskip Let $D := \bigcup \{ D(n): n \geq 1 \}$ be a
$\sigma$-closed-discrete dense set in the underlying LOTS $(X,
\lambda)$. Let $r$ be a monotone countable metacompactness operator
for $(X, \tau)$.  We will show that $R_\tau$ is
$\sigma$-closed-discrete in $(X, \tau)$.   Because $R_\tau
\cap D$ is $\sigma$-closed-discrete in $(X, \lambda)$, it is enough to
show that the set $R' :=
R_\tau - D$ is $\sigma$-closed-discrete in $(X, \lambda)$.

\medskip
For each $p \in R'$, let ${\cal U}(p) := \{ (\leftarrow, p), [p,
\rightarrow)\}$ and find $r({\cal U}(p))$.  Choose any $O(p) \in r({\cal
  U}(p))$ that contains $p$ and note that $O(p) \subseteq [p,
\rightarrow)$ because $r({\cal U}(p))$ refines ${\cal U}(p)$. Then
there is some $y(p) > p$ such that $[p, y(p)) \subseteq O(p)$. Because
$p \in R_\tau$, the set $(p, y(p))$ must be infinite. Let
$N(p)$ be the first integer $k$ such that $(p, y(p)) \cap D(k) \not=
\emptyset$. Decreasing $y(p)$ if necessary, we may assume that $|[p,
y(p)) \cap D(N(p))| = 1$ and we choose the unique $d(p) \in (p, y(p))
  \cap D(N(p))$.

  \medskip For each $k$ and each $d \in D(k)$, let $R(d, k) := \{p \in
  R': N(p) = k, d(p) = d \}$. Note that $R(d,k) \subseteq X - D(k)$
  and it is easy to see that $R(d,k)$ is a subset of the unique convex
  component of $X-D(k)$ that has $d$ as its supremum.

  \medskip It will be enough to show that each set $R(d,k)$ is a
  relatively discrete subspace of $(X, \tau)$, because then Lemma
  \ref{SigmaDiscrete} guarantees that each set $R'(k) := \bigcup \{
  R(d,k): d \in D(k) \}$ is relatively discrete in $(X, \tau)$ and
  therefore (by Lemma \ref{SigmaDiscrete}) is $\sigma$-closed-discrete
  in $(X, \lambda)$. Consequently we will know that the set $R' = \bigcup \{
  R'(k) : k \geq 1 \}$ is also $\sigma$-closed-discrete in $(X,
  \lambda)$ and hence also in $(X, \tau)$, as claimed.

  \medskip Fix any set $R(d_0,k_0)$. If $R(d_0, k_0)$ contains no
  strictly decreasing sequence, then $R(d_0, k_0)$ is well-ordered by
  the given ordering of $X$ and we can write $R(d_0,k_0) = \{
  p(\alpha) : \alpha < \beta \}$.  Because each $p(\alpha) \in R'
  \subseteq R_\tau$, we see that each set $[p(\alpha), p(\alpha + 1))
  \in \tau$ and $[p(\alpha), p(\alpha+1)) \cap R(d_0,k_0) =
  \{p(\alpha)\}$, showing that $R(d_0,k_0)$ is discrete as a subspace
  of $(X, \tau)$, as claimed. Now consider the case where there is
  some strictly decreasing sequence $p(0) > p(1) > \cdots$ in $R(d_0,
  k_0)$.  (We will show that this case cannot occur.) From above,
  $[p(j), d_0) \subseteq O(p(j)) \in r({\cal U}(p(j))$.  Let ${\cal V}
  = \bigcup \{ {\cal U}(p(j)): j \geq 1\}$. Note that ${\cal V}$ is a
  countable open cover of $X$, so $r({\cal V})$ exists. Each ${\cal
    U}(p(j))$ refines $\cal V$ so that $r({\cal U}(p(i))$ refines
  $r({\cal V})$.  Consequently we can choose a set $W(j) \in r({\cal
    V})$ with $$[p(j), d_0) \subseteq O(p(j)) \subseteq W(j).$$
  Note that $p(0) \in W(j)$ for each $j \geq 1$.

\medskip
To complete the proof, we will show that there are infinitely many
distinct sets in the collection $\{W(j): j \geq 1 \}$ and that will
contradict point-finiteness of $r({\cal V})$ at $p(0)$.

\medskip Consider the set $W(1)$.  Because $W(1) \in r({\cal V})$ and
$r({\cal V})$ refines ${\cal V}$, there is some $p(i_1)$ such that
$W(1)$ is contained in some member of ${\cal U}(p(i_1))$, so that
either $W(1) \subseteq (\leftarrow, p(i_1))$ or $W(1) \subseteq
[p(i_1), \rightarrow)$. The first option cannot occur because $W(1)$
contains the non-empty set $(p(1), d_0)$ while $(\leftarrow, p(i_1))$
does not.  Therefore $W(1) \subseteq [p(i_1), \rightarrow)$.  Consider
$j_2 = i_1 + 1$ and the associated set $W(j_2)$.  Because $W(j_2)$
contains $p(j_2)$ while $W(1)$ does not, $W(j_2) \not= W(1)$.
Repeating this argument with $p(j_2)$ and $W(j_2)$ in place of $p(1)$
and $W(1)$ we find $p(i_2)$ with $W(j_2) \subseteq [p(i_2),
\rightarrow)$.  Let $j_3 = i_2 + 1$.  This recursion continues,
producing the required infinite sequence of distinct elements of
$r({\cal V})$ all of which contain $p(0)$, something that is
impossible because $r({\cal V})$ is point-finite.  That completes the
proof that $R_\tau$ is $\sigma$-relatively discrete in $(X,
\tau)$.  

\medskip

The proof that $L_\tau$ is $\sigma$-relatively discrete
is analogous, with ``reverse well-ordering'' in place of
``well-ordering'' when considering the set $L(d_0, k_0)$ (the
analog of $R(d_0, k_0)$ in the above argument). \qed

\medskip

\medskip It is possible that a GO-space $(X, \tau, <)$ has a
$\sigma$-closed-discrete dense set even if the underlying LOTS does
not. For example, if $\delta$ is the discrete topology on the set $X
:= [0, \omega_1)$ with the usual ordering $<$, then the GO-space $(X,
\delta, <)$ has a $\sigma$-closed-discrete dense subset even though
the underlying LOTS, which is the ordinal space $[0, \omega_1)$, does
not. However, a slight modification of the proof of Proposition
\ref{NewFirstHalf} gives:

\medskip
\begin{corollary} \label{FirstHalfCorollary} Suppose the GO-space $(X,
  \tau, <)$ has a $\sigma$-closed-discrete dense subset. If $(X,
  \tau)$ is (countably) monotonically metacompact then $R_\tau \cup L_\tau$ is
  $\sigma$-closed-discrete in $(X, \tau)$. \end{corollary} 

\noindent
Proof: Let $D := \bigcup \{D(k): k \geq 1 \}$ be a
$\sigma$-closed-discrete dense subset of $(X, \tau)$. Then in $(X, \tau)$,
every closed set is a $G_\delta$-set so that every relatively discrete
subset of $(X, \tau)$ is $\sigma$-closed-discrete in
$(X, \tau)$ by Lemma \ref{Perfect}. 

\medskip Use the notation in the proof of Proposition
\ref{NewFirstHalf}. We show that each set $R(d,k)$ is a discrete
subspace of $(X, \tau)$ which makes each set $R(k) := \bigcup \{
R(d,k): d \in D(k) \}$ would also be a discrete subspace of $(X,
\tau)$, and hence $\sigma$-closed-discrete in $(X, \tau)$, so that $R
= \bigcup \{ R(k) : k \geq 1 \} \cup (R \cap D) $ is also be
$\sigma$-closed-discrete in $(X, \tau)$. 

\medskip Fix $(d_0, k_0)$ and consider the set $R(d_0,k_0)$.  As in
the proof of (\ref{NewFirstHalf}), the set $R(d_0, k_0)$ cannot
contain any strictly decreasing sequence, so that it is well-ordered
by the ordering $<$ of $X$ and, just as in (\ref {NewFirstHalf}), must
be relatively discrete, as required. \qed

\medskip Faber's metrization theorem, Theorem 3.1 in \cite{Faber},
will be the key to our next result.  We change some of Faber's notation to
avoid conflicts with the notation used in this paper.

\medskip
\begin{theorem} \label{FabersMetrization} Suppose $(X, \tau, <)$ is a
  GO-space and $Y \subseteq X$. Then the subspace $(Y, \tau_Y)$ is
  metrizable if and only if

\begin{itemize}

\item [{a)}] $(Y, \tau_Y)$ has a dense set $D$ that is the union of
  countably many subsets of $Y$, each being closed-discrete in $(Y,
  \tau_Y)$; and 

\item[{b)}] The sets $\{y \in Y: [y, \rightarrow) \cap Y \in
  \tau_Y \}$ and $\{ y \in Y: (\leftarrow, y] \cap Y \in \tau_Y
\}$ are both $\sigma$-closed-discrete in the subspace $(Y, \tau_Y)$. \qed

\end{itemize}

\end{theorem}

\medskip
\begin{proposition} \label{SecondHalf} Suppose that $(X, \tau, <)$ is a
  GO-space and that the underlying LOTS $(X, \lambda, <)$ has a
  $\sigma$-closed-discrete dense set. If the set $R_\tau \cup L_\tau$ is a
  $\sigma$-discrete subspace of $(X, \tau)$, then
  $(X, \tau)$ is monotonically metacompact. \end{proposition}

\noindent
Proof: Let $D$ be a
$\sigma$-closed-discrete subset of $(X, \lambda)$ that is dense in $(X,
\lambda)$. In the light of Lemma \ref{MakingDiscrete}, the set $R_\tau
\cup L_\tau$ is $\sigma$-closed-discrete in $(X, \lambda)$ and hence
also in $(X, \tau)$. 

\medskip

Let $\mu$ be the topology on $X$ having the collection $$\lambda \cup
\{~ [x, y): x \in R, x<y\} \cup\{~(x, y] : x < y \in L \}$$
as a base.
Then $\mu$ is a GO-topology on $X$ and Faber's metrization theorem
shows that $(X, \mu)$ is metrizable. Also, we see that $\lambda
\subseteq \mu \subseteq \tau$. Now let $S := \{ x \in X: \{x\} \in
\tau \}$.  Then, in the notation of Corollary
\ref{DiscretizingMetric}, $\tau = \mu^S$ showing that $(X, \tau)$ is
monotonically metacompact. \qed

\medskip 
\noindent {\underbar{Proof of Theorems \ref{MainTheorem} and
    Corollary \ref{Metrization}}}: We can combine 
Propositions \ref{NewFirstHalf} and \ref{SecondHalf} to give a proof
of Theorem \ref{MainTheorem} of the Introduction. As noted at the
beginning of the proof of Proposition \ref{SecondHalf}, (c) and (d)
are equivalent. Clearly (a) implies (b) in that theorem, and
Proposition \ref{NewFirstHalf} shows that (b) implies (c).
Proposition \ref{SecondHalf} shows that (c) implies (a).  The proof of
Corollary \ref{Metrization} is similar. Clearly (a) implies (b) in
Corollary \ref{Metrization} and if $(X, \tau)$ is monotonically countably
metacompact, then Corollary \ref{FirstHalfCorollary} shows that
$R_\tau \cup L_\tau$ is $\sigma$-discrete.  In Corollary
\ref{Metrization}, $(X,\tau)$ has a $\sigma$-closed-discrete dense
subset, so that the set $I_\tau$ is also $\sigma$-closed-discrete.
Hence $(X, \tau)$ is metrizable, by Theorem \ref{FabersMetrization}.
\qed

\medskip 
We already proved in Proposition \ref{HereditaryParacompact} that
monotone (countable) metacompactness has certain hereditary
consequences.  A natural question is whether monotone (countable)
metacompactness is itself a hereditary property among GO-spaces.  We
can give an affirmative answer for GO-spaces whose underlying LOTS has
a dense $\sigma$-closed discrete set.  We begin with a lemma.

\medskip
\begin{lemma} \label{IsolatePoints} Suppose $(X, \tau, <)$ is a
  GO-space whose underlying LOTS $(X, \lambda, <)$ has a
  $\sigma$-closed-discrete dense set.  Let $S \subseteq X$ and let
  $\tau^S$ be the topology on $X$ for which $\tau \cup \{ \{x\}: x \in
  S \}$ is a base.  If $(X, \tau, <)$ is monotonically (countably)
  metacompact, then so is the GO-space $(X, \tau^S, <)$. \end{lemma}

\noindent
Proof: As in the Introduction, let $R(\tau) := \{x \in X - I(\tau):
[x, \rightarrow) \in \tau \}$. From Proposition \ref{NewFirstHalf} we
know that $R(\tau)$ is $\sigma$-closed-discrete in $(X, \tau)$ and in
$(X, \lambda)$. Because $\tau \subseteq \tau^S$, we know that
$R(\tau)$ is also $\sigma$-closed-discrete in $(X, \tau^S)$.  In order
to apply Theorem \ref{MainTheorem} to the GO-space $(X, \tau^S, <)$ we
must show that the set $R(\tau^S) := \{ x \in X - I(\tau^S) : [x,
\rightarrow) \in \tau^S \}$ is $\sigma$-closed-discrete in $(X,
\tau^S)$.  But that is automatic because $R(\tau^S) \subseteq
R(\tau)$.  Similarly, the set $L(\tau^S)$ is $\sigma$-closed-discrete
in $(X, \tau^S)$. Now Theorem \ref{MainTheorem} applies to show that
the GO-space $(X, \tau^S, <)$ is monotonically (countably)
metacompact. \qed

\begin{proposition} \label{DiscretizeGO} Suppose $(X, \tau, <)$ is a
  GO-space whose underlying LOTS $(X, \lambda, <)$ has a
  $\sigma$-closed-discrete dense set and suppose that $(X, \tau)$ is
  monotonically (countably) metacompact.  Then for every $Y \subseteq
  X$, the subspace $(Y, \tau_Y)$ of $(X, \tau)$ is also monotonically
  (countably) metacompact. \end{proposition}

\noindent
Proof: Let $Y \subseteq X$.  Let $S = X - Y$ and create the topology
$\tau^S$ as in Corollary \ref{DiscretizingMetric}.  By Lemma 
\ref{IsolatePoints} we know that $(X, \tau^S)$
is monotonically (countably) metacompact.  Note that $(\tau^S)_Y =
\tau_Y$, i.e., that $(Y, \tau_Y)$ is a subspace of $(X, \tau^S)$.  In
fact, $(Y, \tau_Y)$ is a {\em{closed}} subspace of the monotonically
(countably) metacompact space $(X, \tau^S)$ and therefore inherits
monotone (countable) metacompactness. \qed

\medskip
\section{Examples and questions}

\medskip

\medskip Suppose $(X,<)$ is a linearly ordered set and $Y \subset X$.
We say that a set $S \subset Y$ is {\em{relatively convex in $Y$}} if
a point $b$ of $Y$ has $b \in S$ whenever $a \leq b \leq c$ for points
$a,c \in S$.  For any subset $T \subseteq Y$ we let $C(T) = \bigcup \{
[u,v] : u \leq v, u,v \in T \}$ and we refer to $C(T)$ as the {\em{convex
    hull of  $T$ in $X$}}.  Note that $C(T) \cap Y = T$ provided $T$ is
a relatively convex subset of $Y$.

\begin{example} \label{MonotmetaLOTS} There is a first-countable
  monotonically metacompact LOTS that is not metrizable. \end{example}

\noindent
Proof: Let $X := (\IR \times \{0\} ) \cup (\IP
\times \IZ)$ with the lexicographic ordering and
the open interval topology $\lambda$ of that
ordering. This LOTS $(X, \lambda)$ contains the
Michael line and is therefore non-metrizable.

\medskip Let $Y = \IR \times \{0\}$ and let $\mu$ be the subspace
topology that $Y$ inherits from $(X, \lambda)$.  Then $(Y, \mu, <)$ is
the Michael line so that Proposition \ref{SecondHalf} gives a monotone
metacompactness operator $r_Y$ for $(Y, \mu)$. We may assume that
$r_Y$ always produces collections of sets that are relatively convex
in $Y$. Let $\cal U$ be any open cover of $X$. Without loss of
generality, we may assume that members of $\cal U$ are convex open
subsets of $X$. Let ${\cal U}_Y := \{ U \cap Y: U \in {\cal U} \}$.
Find $r_Y({\cal U}_Y)$. For any $S \in r_Y({\cal U}_Y)$ let $C(S)$ be
the convex hull of $S$ in $X$. Because $S$ cannot contain a rational
endpoint of itself, it is easy to check
that each $C(S)$ is open in $X$. The
collection $r_1({\cal U}) := \{C(S): S \in r_Y({\cal U}_Y) \}$ refines
$\cal U$, covers $Y$, and is point-finite in $X$. To complete the
proof, let $r({\cal U}) := r_1({\cal U}) \cup \{ \{(x,n)\} \in X: n
\not= 0 \}$. \qed

\medskip
With a little more care, we can construct a non-metrizable LOTS that
is monotonically metacompact and monotonically
Lindel\"{o}f. We thank 
Dennis Burke for pointing out the next example.

\medskip
\begin{example} There is a non-metrizable LOTS ~$X$ having a monotone
  metacompactness operator $R$ with the additional property that for
  each open cover ${\cal U}$ of $X,~ R({\cal U})$ is countable.  Hence
  $R$ is also a monotone Lindel\"{o}f operator in the sense of
  \cite{BennettLutzerMatveev}. \end{example}

\noindent
Proof: Let $B \subseteq \IR$ be a Bernstein set, i.e., a subset of
$\IR$ such that neither $B$ nor $C := \IR - B$ contains an uncountable
compact subset of $\IR$.  Let $X = (\IR \times \{0\}) \cup (C \times
\IZ)$ have the open interval topology $\lambda$ of the lexicographic
order. Let $Y = \IR \times \{0\}$ and let $\mu$ be the usual open
interval topology on the set $Y$. Note that $\mu$ is {\em{not}} the
subspace topology that $Y$ inherits from $(X, \lambda)$.

\medskip Let ${\cal U}$ be an open cover of $X$ by convex sets. Let
$a({\cal U}) := \{ \Int_\mu( U \cap Y) : U \in {\cal U}\}$.  Then
$a({\cal U})$ is a collection of open sets in the metric space $(Y,
\mu)$ that covers $B$.  According to Proposition \ref{MetaMoore} there
is a monotone operator $b$ such that $b({\cal U})$ refines $a({\cal
  U})$, has $\bigcup b({\cal U}) = \bigcup a({\cal U})$, and is a
point-finite collection of relatively convex open intervals in $Y$.
Because $(Y, \mu) $ is separable, $b({\cal U})$ must be countable.
Because $b({\cal U})$ covers the Bernstein set $B \times \{0\}$, the
set $Y - \bigcup b({\cal U})$ is countable. For any $S \in b({\cal
  U})$ let $C(S)$ be the convex hull of $S$ in $X$. Let $$R({\cal U})
:= \{ C(S): S \in b({\cal U}) \} \cup \{ \{(x,n)\}: (x,n) \in X ~
{\rm{and}}~ (x,0) \in Y - \bigcup b({\cal U}) \}.$$  Then $R({\cal U})$
is the required countable, open, point-finite refinement of $\cal U$. \qed

\medskip
Although most of our results use the existence of
$\sigma$-closed-discrete dense sets, they can sometimes be applied in
more general contexts.

\begin{example} The lexicographic square $X = [0,1] \times [0,1]$ with
  the open interval topology of the lexicographic order is not
  monotonically (countably) metacompact. \end{example}

\noindent
Proof: The space $X$ contains the Alexandroff double interval $Y :=
[0,1] \times \{0, 1\}$ as a closed subspace. Proposition
\ref{SecondHalf} shows that $Y$ is not monotonically countably
metacompact. Hence neither is $X$. Alternatively, apply Proposition
\ref{CompactMonoMeta}.\qed

\medskip

Our results characterize monotone metacompactness in GO-spaces whose
underlying LOTS has a $\sigma$-closed-discrete dense set.  Souslin
lines are historically important examples
of LOTS that are perfect but do not have any $\sigma$-closed-discrete
dense subsets.  Whether or not Souslin lines exist is undecidable in
ZFC \cite{Todorcevic}. As our next two examples show, in any model of
ZFC that contains Souslin lines, some Souslin lines will be monotonically
metacompact, and others will not be.

\begin{example} \label {SouslinYes} If there is a
  Souslin line, then there is a Souslin line
  that is monotonically metacompact. \end{example}

\noindent
Proof: If there is a Souslin line, then there is a Souslin tree $(T,
\preceq)$, i.e., a tree with uncountably many levels in which all each
levels and all anti-chains are countable.  Then there is a Souslin
tree $T$ with the properties that each $t \in T$ has $\omega$-many
immediate successors \cite{Todorcevic}. As in \cite{Todorcevic}, a
{\em{node of the tree}} is the set of all members of the tree with exactly
the same set of predecessors.  Each node is an anti-chain, so each
node is countable. Choose an ordering for each
node that makes it look like the set of all integers.

\medskip
Let $B(T)$ be the branch space of $T$. Then $B(T)$
is linearly ordered by the ``first-difference
ordering''. Put the open interval topology of that
linear ordering on $B(T)$. Because of the special 
node-orderings chosen above, for each $t \in T$
the set $[t] := \{ b \in B(T): t \in b \}$ is a
convex open set in the branch space and the set of
all possible $[t]$ is a basis for $B(T)$.

\medskip
Let ${\cal U}$ be any open covering of $B(T)$.
Given $b \in U \in {\cal U}$ there is some $t \in
T$ with $b \in [t] \subseteq U$. We will say that $t$
is ${\cal U}$-minimal if $[t]$ is contained in
some member of $\cal U$ and if no predecessor of
$t$ in $T$ has this property.  Let $r({\cal U}) :=
\{ [t] : t \in T ~{\rm{and}}~ t ~{\rm{is}}~ {\cal
  U}{\rm{-}}minimal\}$. Then $r({\cal U})$ is an open
cover of $B(T)$ that refines ${\cal U}$.  Note that
if $[t_1], [t_2] \in r({\cal U})$ and $b \in B(T)$
  has $b \in [t_1] \cap [t_2]$, then $t_1, t_2 \in
  b$.  Because any two members of the branch $b$
  are comparable in $T$, we have either $t_1 \prec
  t_2$ or vice-versa, and both are impossible by
  ${\cal U}$-minimality of $t_1$ and $t_2$. \qed

\begin{remark} Note that the argument in the previous example shows
  that any non-Archimedean space is monotonically metacompact.\end{remark}

\begin{example} \label{SouslinNo} If there is a
  Souslin line, then there is a Souslin line that
  is not monotonically metacompact. \end{example}

\noindent
Proof: We give two examples, one connected and one totally
disconnected. If there is a Souslin line, then there is a
compact connected Souslin line $S_1$ and a compact totally
disconnected Souslin line $S_2 = S_1 \times \{0,1\}$ with the
lexicographic ordering. Because $S_1$ and $S_2$ are compact LOTS that
are not metrizable, Proposition 
\ref{CompactMonoMeta} shows that neither is monotonically
metacompact. \qed

\begin{question} Characterize monotone (countable) metacompactness in
  GO-spaces, without making assumptions about the existence of special
  dense sets (as in Proposition \ref{NewFirstHalf} and Corollary
  \ref{Metrization}). \end{question}

\begin{question} Must the GO-space $(X, \tau,<)$ be monotonically
  (countably) metacompact if its subspace $Y := X - I_\tau$ is
  monotonically (countably) metacompact? \end{question}

\begin{question} If $(X, \tau,<)$ is a monotonically (countably)
  metacompact GO-space and $Y \subseteq X$, must the subspace $(Y,
  \tau_Y)$ be monotonically (countably) metacompact? In other words,
  is monotone (countable) metacompactness a hereditary property among
  GO-spaces? By Proposition \ref{AlmostHereditary}, the answer is
  ``Yes'' in case the underlying LOTS of $(X, \tau, <)$ has a
  $\sigma$-closed-discrete dense set, but in general this question
  remains open. \end{question}

\begin{question} If $(X, \tau,<)$ is a monotonically (countably)
  metacompact GO-space and $S \subseteq X$, is the GO-space $(X,
  \tau^S, <)$ also monotonically (countably) metacompact? (See Lemma
  \ref{IsolatePoints} for the definition of $\tau^S$.) Note that the
  proof of Proposition \ref{AlmostHereditary} shows that an
  affirmative answer to this question would give an affirmative answer
  to the previous question. \end{question}

\begin{question} Suppose $X$ is a compact Hausdorff space that is
monotonically metacompact.  Is $X$ metrizable? (Compare Gruenhage's
theorem that a monotonically compact Hausdorff space is metrizable
\cite{Gruenhage}.) \end{question}

In Proposition \ref{MetaMoore} we showed that
among Moore spaces, metacompactness and monotone
metacompactness are equivalent properties. This
suggests investigating the role of monotone
metacompactness in other generalized metric spaces.

\begin{question} Must a metacompact quasi-developable space $X$ be
  monotonically metacompact? What if $X$ is hereditarily metacompact
  (so that each level of the quasi-development may be assumed to be
  point-finite)? What if $X$ has a $\sigma$-disjoint base?
\end{question} 

\begin{question} Which stratifiable spaces \cite {Borges} are
  monotonically metacompact? Two particularly interesting examples of
  stratifiable spaces are due to McAuley and Ceder (see
  \cite{AmonoMizokami} for a description of these spaces).  Are these
  spaces monotonically metacompact?
\end{question}

\begin{question} \label{RefereeQuestion} The referee suggested that we
  ask which of our results can be proved for the class of
  monotonically normal spaces (which is wider than the class of
  GO-spaces).  We have already remarked that a monotonically normal
  space that is monotonically countably metacompact must be
  paracompact.  We do not know whether such a space is
  {\em{hereditarily}} paracompact. In addition, we do not know whether
  a compact monotonically normal space must be metrizable whenever it
  is monotonically countably metacompact.\end{question}

\end{document}